\documentclass[12pt,a4paper]{article}
\usepackage{amsfonts, amsmath, amscd, amssymb}
\usepackage[T2A]{fontenc}
\usepackage[cp1251]{inputenc}
\usepackage[english]{babel}

\newcommand{\pr}{\par {\bfseries\itshape Proof.} }   
\newcounter{df}
\newcounter{sen}
\newenvironment{sen}{\par
\refstepcounter{sen}
{\bf Proposition \arabic{sen}.} }{}
\newcounter{rem}
\newenvironment{rem}{\par
\refstepcounter{rem}
{\bfseries\itshape Remark \arabic{rem}.} }{}
\newcounter{exa}
\newcounter{teo}
\newenvironment{teo}{\par
\refstepcounter{teo}
{\bf Theorem \arabic{teo}.} \it }{}
\newcounter{cor}
\newenvironment{cor}{\par
\refstepcounter{cor}
{\bfseries\itshape Corollary \arabic{cor}.} \it }{}
\newcounter{st}
\newcounter{lem}
\newenvironment{lem}{\par
\refstepcounter{lem}
{\bf Lemma \arabic{lem}.} \it }{}
\makeatletter                                                   %
\renewcommand{\section}{\@startsection{section}{1}
{\parindent}{3.5ex plus 1ex minus 0.2ex}{2.3ex plus 0.2ex}{\bf}}%
\makeatother                                                    %
\makeatletter                                                   %
\@addtoreset{equation}{section}                                 %
\makeatother                                                    %
\makeatletter                                                   %
\@addtoreset{exa}{section}                                      %
\makeatother                                                    %
\makeatletter                                                   %
\@addtoreset{teo}{section}                                      %
\makeatother                                                    %
\makeatletter                                                   %
\@addtoreset{df}{section}                                       %
\makeatother                                                    %
\makeatletter                                                   %
\@addtoreset{cor}{section}                                      %
\makeatother                                                    %
\makeatletter                                                   %
\@addtoreset{lem}{section}                                      %
\makeatother                                                    %
\makeatletter                                                   %
\@addtoreset{rem}{section}                                      %
\makeatother                                                    %
\makeatletter                                                   %
\@addtoreset{sen}{section}                                      %
\makeatother                                                    %
\makeatletter                                                   %
\@addtoreset{st}{section}                                       %
\makeatother                                                    %

\author{S.S. GABRIYELYAN\footnote{The author was partially supported
 by Israel Ministry of Immigrant Absorption and ISF grant 888/07}}
\title{ABSOLUTE CONTINUITY AND SINGULARITY OF TWO PROBABILITY
MEASURES ON A FILTERED SPACE}
\date{}

\begin{document}

\makeatletter
\renewcommand{\@makefnmark}{}
\renewcommand{\@makefntext}[1]{\parindent=1em #1}
\makeatother

\maketitle

\vspace{-\baselineskip}

\begin{abstract}
Let $\mu$ and $\nu$ be fixed probability measures on a filtered
space $(\Omega , {\cal F}, ({\cal F}_t )_{t\in {\bf R}^{+} } )$.
Denote by $\mu_T $ and $\nu_T $ (respectively, $\mu_{T-} $ and
$\nu_{T-} $) the restrictions of the measures $\mu$ and $\nu$ on
${\cal F}_T $ (respectively, on ${\cal F}_{T-} $) for a stopping
time $T$. We  find the Hahn decomposition of $\mu_T $ and $\nu_T
$ using the Hahn decomposition of the measures $\mu$, $\nu$, and the
Hellinger process $h_t$ in the strict sense of order
$\frac{1}{2}$. The norm of the absolutely continuous component of
$\mu_{T-} $ with respect to $\nu_{T-} $ is computed in terms of density processes
and Hellinger integrals.
\end{abstract}

{\bf Introduction.}
Let $M^{+}(\Omega )$ be the set of all nonnegative finite
measures on a measurable space $\Omega $. The norm (=the total variation) of a measure $\mu$ is denoted by $\|\mu\|$. A measure $\mu\in M^{+}(\Omega )$ is
called a probability measure if $\mu (\Omega ) =1$. For $\mu ,\nu \in
M^{+} (\Omega )$ we write $\mu\ll\nu $ (resp.
$\mu\perp\nu $) if $\mu $ is absolutely continuous (resp. singular)
 with respect  to $\nu $. Mutual absolute continuity (equivalence) of
$\mu$ and $\nu $ we denote by $\mu\sim\nu $. If $\mu = \mu_1 +
\mu_2 $, with $\mu_1 \perp \mu_2 $, then $\mu_1 $ and $\mu_2 $
are called parts of $\mu $.
Let $\mu ,\nu \in M^{+}(\Omega )$. We can write  the Lebesgue decomposition of the measures $\mu$ and $\nu $  with respect
to each other in the following form
\[
 \mu = \mu^1 + \mu^2 , \nu = \nu^1 + \nu^2 , \mbox{ with }
\mu^1 \sim \nu^1 , \mu^2 \perp \nu , \nu^2 \perp \mu .
\]
  Denote the derivative of $\mu$ with respect to
$\nu $ by $\frac{d\mu }{d\nu }$. Then
\[
 \frac{d\mu }{d\nu } = \frac{d\mu^1 }{d\nu^1 } ,\enspace \nu^1
-\mbox{a.s.}; \enspace \mbox{ and } \enspace \frac{d\mu }{d\nu }
= 0 , \enspace (\nu^2 + \mu^2 )-\mbox{a.s.}
\]

Let  $\Omega =(\Omega , {\cal F}, \mathbf{F} =({\cal F}_t )_{t\in {\bf R}^{+} } )$ be  a filtered space  with a right continuous filtration and ${\cal F} = \vee_t {\cal F}_t $,  $T$ be a stopping time and  $\mu$ a probability measure on $\Omega$. We denote by
$\mu_{T-} $ and $\mu_T $ the restrictions of $\mu $  on ${\cal F}_{T-} $ and   ${\cal F}_T $ respectively. For $E\in {\cal F}_{T}$ the stopping time $T_E$ is defined as follows: $T_E (\omega) = T(\omega)$, if $\omega\in E$, and $T_E (\omega) = +\infty$, if $\omega\not\in E$.

Fix two  probability measures $\mu$ and $\nu$ on $\Omega$ and set  ${\rm Q} =\frac{1}{2} (\mu + \nu )$. Denote by $z$ and $z' $ the density processes of $\mu $ and $\nu $ with respect to ${\rm Q}$ which are right-continuous and admit left-hand limits.  We define the process $z_{-}$ as follows: $z_{0-} =z_0$ and $z_{t-} =\lim_{s\to t-0} z_s$. Analogously, the process $z'_{-}$ is defined.
We shall write $\mu \stackrel{loc}{\ll } \nu $ if  $\mu $ is locally absolutely continuous with respect to $\nu $, i.e., if $ \mu_t \ll \nu_t ,  \forall t$.
The set $J=\{ \beta : \beta \mbox{ is a part of } \mu  \mbox{ such that } \beta \stackrel{loc}{\ll } \nu \}$ is a partially ordered set with the following natural binary relation: for $ \beta_1 ,\beta_2 \in J$,  $\beta_1 \leq \beta_2$ iff $\beta_1$ is a part of $\beta_2$. By the Zorn lemma, $J$ contains the maximal element $\overline{\mu}$ that  is called the {\it locally absolutely continuous part} of $\mu $ with respect to $\nu $.  The part ${\tilde {\tilde\mu }} = \mu - \overline{\mu}$ of
$\mu $ we call the {\it asymptotic singular part} of $\mu $ with respect
to $\nu $ (justification of the title ''asymptotic singular part'' is contained in Proposition  \ref{l4}). We shall write $\mu \stackrel{as}{\perp } \nu $ if $\overline{\mu}  =0$.

Let $S_n = \inf ( t : \ z_t <\frac{1}{n} $ or $z'_t <\frac{1}{n} ) $ (where $\inf \varnothing :=\infty$). The stopping time
$S$ is the first moment when either $z$ or $z' $ vanishes,
\[
 S= \inf ( t : \  z_t =0 \mbox{ or } z'_t =0 ) .
\]

The process $ Y(\alpha ) = z^{\alpha } {z'}^{1-\alpha } $, where
$\alpha \in (0;1)$, is a ${\rm Q}$-supermartingale of class
$(D)$. If $\alpha =0.5$ we shall write $Y_t =
\sqrt{z_t z'_t }$. The number $ H(\alpha ; \mu ,\nu ) = {\bf
E}_{\rm Q} [ Y(\alpha )_\infty ] $ is
called the Hellinger integral of the order $\alpha$. Let $Y= M-A$ be the Doob-Meyer decomposition of $Y$. By \cite[IV.1.18]{JS}, there exists a predictable increasing $\overline{\mathbb{R}}_+$-valued process  $h_t$, unique up to ${\rm Q}$-indistinguishability, such that $h_0 =0$ and the following two conditions hold
\begin{equation} \label{f1}
A= Y_{-} \bullet h \; , \quad h=\left( \frac{1}{Y_{-} }
1_{\Gamma'' } \right) \bullet A , \mbox{ where } \Gamma'' =\left( \{ z_{-} >0\} \cap \{ z'_{-} >0\}\right) \cup [0].
\end{equation}
The process $h$ is called  the Hellinger process in
the strict sense,  of order $\frac{1}{2}$, between $\mu$ and $\nu$.
 The Hellinger process $h(1)$ of order $1$ is defined as the ${\rm
Q}$-compensator of the following process (see \cite[ IV.1.53]{JS}, where
$0/0=0$)
\begin{equation} \label{f2}
A^1 = \frac{z_S }{z_{S-}} 1_{\{ 0<S<\infty , \  z'_S =0 < z'_{S-}
\} } 1_{ \mbox{\boldmath$[$} S, \infty \mbox{\boldmath$[$} } .
\end{equation}

A stopping time $T$ is called the {\it stopping time of a process}
$X$ if: 1) $X=X^T$, 2) if $X= X^U$ for a stopping time $U$, then $T\leqslant U$, ${\rm Q}$-a.s. Taking $\textrm{essinf}$ (${\rm Q}$-a.s.) of the stopping times $U$
with the property $X = X^U$ we can see that for any  right continuous process there exists
its stopping time. Importance of this notion is
demonstrated in Theorem \ref{th2}.

Let $X$ be a process and $T$ be a stopping time. Taking into
account the evident physical interpretation: the process $X^{T-} =
X 1_{[0; T[ } $ is called {\it the process $X$ interrupted at the
moment $T$}.

If sets $E_1$ and $E_2$ are mutually disjoint, their union is denoted by $E_1 \sqcup E_2$.
A decomposition $\Omega = E\sqcup E^c $, where $E^c
=\Omega\setminus E$, is called the Hahn decomposition of measures
$\mu $ and $\nu $ if: $1)\ \mu\sim\nu $ on the set $E$; $\; 2)\
\mu\perp\nu $ on the set $E^c$. Note that the Hahn decomposition is unique up to ${\rm Q}$-negligible sets.

The question of absolute continuity or
singularity of two probability measures has been investigated a
long time ago, both for its theoretical interest and for its
applications to mathematical statistics, financial mathematics,
ergodic theory and others. S.~Kakutani in 1948 \cite{Kak}, was the
first to solve this problem in the case of two measures having an
infinite product form. Yu.~M.~Kabanov, R.~Sh.~Liptser, A.~N.~Shiryaev
\cite{KLS, KLS2}(see also \cite[ \S 6, ch. 7]{Sh})
generalized this result for measures on the $\sigma $-algebra
${\cal B} $ which is generated by an increasing sequence of
$\sigma $-algebras ${\cal B}_n $ (under the condition of their
local absolute continuity). A.~R.~Darwich \cite{Dar} extended
Theorem 4 of Yu.~M.~Kabanov et al. \cite{KLS}.
The following question, which has been
considered by several authors,  is the main theme of the chapter
IV of the book \cite{JS}:
\begin{enumerate} \item[]
{\bf Problem 1.} {\it Under which conditions can we assert that
$\mu_T \ll\nu_T$ or $\mu_T\perp\nu_T$}?
\end{enumerate}

This problem can be attacked via "Hellinger integrals"
and  "Hellinger processes". However, a situation may naturally
occur, where the two measures are neither (locally) absolutely
continuous nor singular. W.~Schachermayer and W.~Schachinger
\cite{SS} have raised the more general question:
\begin{enumerate} \item[]
{\bf Problem 2.} {\it Can we find the Hahn decomposition of $\mu_T $
and $\nu_T $?}
\end{enumerate}
In \cite{JS} and \cite{SS} the authors have looked for the answers
to these questions using the values of the Hellinger processes of
different orders at time $T$ (i.e., in "predictable" terms) and the difficulty to find them arises
from the fact that $z_t$ may jump to zero. Unfortunately, there is no hope to obtain a complete answer to Problem 2, as the following simple example, which has been constructed in \cite{SS}, shows. There exists a filtered space equipped with two probability measures $\mu$ and $\nu$ with the following property: there is no $[0,\infty]$-valued predictable process $H$ such that, for every stopping time $T$, $\nu_T \perp \mu_T \Leftrightarrow  \nu(H_T =\infty) =1$. Such an example may be constructed since in general $\mathcal{F}_{T-} \not= \mathcal{F}_T$ and hence $z_t$ may suddenly jump to zero. Nevertheless, W.~Schachermayer and W.~Schachinger \cite{SS} have obtained the following interesting positive result:

{\bf Theorem A}. {\it Under the assumption that $\mu_0 \sim \nu_0$ we have, for every stopping time $T$, }
\begin{equation} \label{e01}
\{ S\leqslant T, z_{S-} =0 \} = \{ h_T =\infty \}, \quad \nu\mbox{-a.s.}
\end{equation}

Clearly the stopping time $S$ plays an important role.
It is easy to give a simple answer to Problem 2  if we know $S$
and the set $B= \{ 0< z_{\infty} <2\}$ on which $\mu \sim \nu$.
A.~S.~Cherny and M.~A.~Urusov \cite{Che} added a point $\delta$ to
$[0;\infty]$ in such a way that $\delta >\infty$ and considered
the  separating time $\widetilde{S}$ for $\mu$ and $\nu$:
\[
\widetilde{S} (\omega)=S(\omega) \mbox{ if } \omega\in B^c \mbox{
and } \widetilde{S} (\omega)=\delta \mbox{ if } \omega\in B.
\]
The following theorem is proved in \cite{Che}:

{\bf Theorem B}. {\it For any stopping time $T$ we have}
\[
\mu_T \sim \nu_T \mbox{ {\it on the set} } \{ T< \widetilde{S} \}
\mbox{ and } \mu_T \perp \nu_T \mbox{ {\it on the set} } \{ T\geqslant
\widetilde{S} \} .
\]
In \cite{Che2} the authors computed of $\widetilde{S}$ in many
important cases.

If the stopping time $S$ and the process $h$ are known, we give the following answer to
Problem 2:
\begin{teo} \label{th1}
{\it Let $T$ be a stopping time. Set}
\[
E\ =\left( \{ T<S \} \cup \{ T=S ,\ T=\infty \} \right) \cap \{
h_T < \infty \}
\]
\[
 E^c =\left( \{ S<T \} \cup \{ S\leqslant T ,\ T<\infty \}
\right) \cup \{ h_T = \infty \} .
\]
{\it Then: $\mu_T \sim\nu_T $ on the set $E$, and $\mu_T \perp\nu_T
$ on the set $E^c $.}
\end{teo}
\\ In particular, if $\mu \stackrel{loc}{\ll } \nu $, then
$S\equiv\infty $ and Corollary IV.2.8 of \cite{JS} follows from
Theorem \ref{th1}.


Theorem \ref{th1} shifts Problem 2 into the next problem:
\begin{enumerate} \item[]
{\bf Problem 3.} {\it Find the stopping time S}.
\end{enumerate}
If  $h$ and the Hahn decomposition of the measures $\mu$ and $\nu$ are known, we propose the following solution of Problem 3:

\newpage
\begin{teo} \label{th2}
{\it
\begin{enumerate}
\item The stopping time $H$ of the process $h$ exists and it coincides with the stopping times of the processes $A, M, Y, z$ and $z' $. Moreover, (${\rm Q}$-a.s.)
\begin{equation} \label{f3}
H\leqslant S  \; \mbox{ and } \quad \{ H<S\} = \{ 0< z_H < 2, H<\infty \}
\subseteq \{ S=\infty \} .
\end{equation}
\item $\mu\sim\nu $ on the set $\{ H<S\} $, and $S= H_{\{ H=S\} }
$.
\item There exists a version $B\sqcup B^c $ of the Hahn decomposition of the
measures $\mu $ and $\nu $ with $B^c \in {\cal F}_H $, where $\mu\sim\nu $ on $B$, such that (${\rm Q}$-a.s.)
\[
S= H_{B^c }.
\]
\end{enumerate} }
\end{teo}

\begin{rem}
Set $S^0 =S_{ \cup_n \{ S_n =S\} } $. By Definition IV.1.24 of \cite{JS}
and Theorem A,  any Hellinger process of order
$\frac{1}{2} $ is equal to $h(\frac{1}{2} ; \mu , \nu ) = h + A'
1_{ \mbox{\boldmath$]$} S^0 , \infty \mbox{\boldmath$[$} } $,
where $A'$ is a predictable increasing process. Then the stopping
time $H'$ of the process $h' =h + t 1_{ \mbox{\boldmath$]$} S^0 ,
\infty \mbox{\boldmath$[$} } $ is equal to $H_{ \{ H<S\} \cup \{
h_H = \infty \} }$. This is the greatest stopping time of
Hellinger processes of order $\frac{1}{2}$ while $H$ is the
smallest one. Obviously,  $S= H_{ \{ H<H' \} \cup \{ h_H
= \infty \} }$. So it is also important to find another Hellinger processes which have the
greatest stopping time.
\end{rem}

Theorem \ref{th2} shows the importance of the stopping time $H$
of the process $h$. We show that the knowledge of the stopping time
$H_1$ of the process $h(1)$ does not determine $S$ either.
\begin{sen} \label{p11}
{\it Set $N_1 = \{ 0<S<\infty ,\ z'_S =0< z'_{S-} \} $ and
\[
H_1 = \mathrm{essinf} \{ W : \ W \mbox{ is a stopping time such that } N_1 \in \mathcal{F}_W \mbox{ and } W_{N_1 } = S_{N_1 }, \; {\rm Q}-\mbox{a.s.} \}.
\]
Then $H_1 $ is the
 stopping time  of $h(1)$.}
\end{sen}

If $T=\infty$, then ${\cal F}_{T} ={\cal F}_{T-} $. Thus, it is natural to consider the counterparts of Problems 1 and 2 for $\mu_{T-} $ and $\nu_{T-} $. In what follows we consider the next question:
\begin{enumerate} \item[]
{\bf Problem 4.} {\it Find the norm of the absolutely continuous
component of $\mu_{T-} $ with respect  to $\nu_{T-} $. }
\end{enumerate}

In the following theorem we give the solution of Problem 4 (in terms of density
processes and Hellinger integrals). Note that for this theorem it is enough to know only the density
processes $z^{T-}$ and $z'^{T-}$ interrupted at the moment $T$;
$z_0 , z'_0 $ and the system ${\cal L}= \{ {\cal F}_0 \mbox{ and
} A\cap \{ t<T \} , A\in {\cal F}_t \} $ which generates ${\cal
F}_{T-} $.

\begin{teo} \label{th3}
{\it Let probability measures $\mu , \nu $ and ${\rm P}$ on a filtered space
$(\Omega , \mathcal{F}, \mathbf{F} )$ be such that $\mu
\stackrel{loc}{\ll } {\rm P} , \nu \stackrel{loc}{\ll } {\rm P}$ and
let $z$ and $z' $ be the density processes of $\mu $ and $\nu $
 with respect  to ${\rm P}$ respectively. If $T$ is a stopping time, then$^1$\footnotemark\footnotetext{${}^1{\underline\lim }_{n\to\infty \atop
x\to 1-0} f(x,n) = \sup_{n\in \mathbb{N} \atop \alpha \in (0,1)} \inf_{m\geqslant n \atop x\in [\alpha , 1)} f(x,n)$.}
\[
 \| (\mu_{T-} )_a \| = {\underline\lim }_{n\to\infty \atop
\alpha\to 1-0} \left\{ \int z_0^{\alpha } {z'}_0^{1-\alpha }
1_{\{ T=0\} } d{\rm P}_0 + \int z_n^{\alpha } {z'}_n^{1-\alpha }
1_{\{ n<T\} } d{\rm P}_n + \right.
\]
\begin{equation} \label{f5}
\left. \sum_{k=1}^{n2^n } \int \left[ {\bf E}_{\rm P} \left[
z_{\frac{k}{2^n } } 1_{ \left\{ \frac{k-1}{2^n } < T\leqslant
\frac{k}{2^n } \right\} } \big | {\cal F}_{\frac{k-1}{2^n } }
\right] \right]^{\alpha } \left[ {\bf E}_{\rm P} \left[ z'_{\frac{k}{2^n }
} 1_{ \left\{ \frac{k-1}{2^n } < T\leqslant \frac{k}{2^n } \right\} }
\big | {\cal F}_{\frac{k-1}{2^n } } \right] \right]^{1-\alpha }
d{\rm P}_{\frac{k-1}{2^n } } \right\} ,
\end{equation}
where $(\mu_{T-} )_a $ is the absolutely continuous part of
$\mu_{T-} $  with respect  to $\nu_{T-} $. }
\end{teo}

Our proof of Theorem \ref{th3} follows in three steps:

{\bf Step 1.} {\it We prove Theorem {\rm \ref{th3}} assuming that
$T\equiv\infty$ and the time-set is $\mathbb{N}$.}

More precisely, setting
\[
H(\alpha ; \mu_n ,\nu_n ) =\int \left( \frac{d\mu_n }{d\nu_n } \right)^{\alpha }
d\nu_n ,
\]
we prove the following:
\begin{sen} \label{p0}
{\it Let $T\equiv\infty $ and $\alpha \in (0;1)$.  Then
\[
 \| \mu_a \| = {\underline\lim }_{n\to\infty \atop \alpha\to 1-0} H(\alpha ; \mu_n ,\nu_n ),
\]
 where $\mu_a $ is the absolutely continuous part of $\mu $
 with respect  to $\nu $.}
\end{sen}

To prove Proposition \ref{p0} we essentially use Jessen's theorem \cite[Theorem 5.2.26]{Str} which allows us to find even the density of $\mu$  with respect  to $\nu$:

{\bf Theorem C}. (Jessen) {\it Let $\mu_a $ be the absolutely continuous part of $\mu $
 with respect  to $\nu $. Then
\[
 \frac{d\mu_a }{d\nu } =\lim_{n\to\infty } \frac{d(\mu_n)_a }{d\nu_n
}, \quad \nu\mbox{-a.s.}
\] }

Before proceeding to step 2, let us note  a few immediately corollaries of Proposition \ref{p0} and the Jessen theorem.

For a predictable stopping time we can
compute the norm of the absolutely
continuous part of $\mu_{T-}$  with respect  to $\nu_{T-}$ as follows:
\begin{cor} \label{c0}
{\it Let $\mu $ and $ \nu $ be two probability measures on
 a filtered space $(\Omega , \mathcal{F}, \mathbf{F} )$. If a stopping time $T$ is
predictable and a sequence $\{ V_n \} $ is an announcing sequence
for $T$, then
\[
 \| (\mu_{T-} )_a \| = {\underline\lim
}_{n\to\infty \atop \alpha\to 1-0} H(\alpha ; \mu_{V_n } ,
\nu_{V_n } )\; \mbox{ and } \; \frac{ d(\mu_{T-} )_a }{d\nu_{T-} } =
\lim_{n\to\infty } \frac{ d(\mu_{V_n } )_a }{d\nu_{V_n } }, \quad
\nu_{T-}\mbox{-a.s.},
\]
 where $(\mu_{T-} )_a $ is the absolutely
continuous part of $\mu_{T-} $  with respect  to $\nu_{T-} $. }
\end{cor}

Note that a similar result for $\mu_{T}$ and $\nu_{T}$ can not be established, since the filtration is not left continuous in the general case (i.e., if $\mathcal{F}_T \not= \mathcal{F}_{T-}$, then absolute continuity and singularity at the moment $T$ are not determined by the preceding events). We can obtain only the following analog of Corollary \ref{c0}:
\begin{cor} \label{c1}
{\it Let $\mu $ and $ \nu $ be probability measures on a filtered space $(\Omega ,
\mathcal{F}, \mathbf{F} )$ and let a nondecreasing sequence
$\{ V_n \} $ of stopping times be such that $\lim_n V_n =\infty
$. If $T$ is a stopping time, then
\[
 \| (\mu_{T} )_a \| =
{\underline\lim }_{n\to\infty \atop \alpha\to 1-0} H(\alpha ;
\mu_{T\wedge V_n } ,  \nu_{T\wedge V_n } )\; \mbox{ and } \; \frac{
d(\mu_T )_a }{d\nu_T } = \lim_{n\to\infty } \frac{ d(\mu_{T\wedge
V_n })_a }{d\nu_{T\wedge V_n } }, \; \nu_{T}\mbox{-a.s.},
\]
where $(\mu_{T} )_a $ is the absolutely continuous part of
$\mu_{T} $  with respect  to $\nu_{T} $. }
\end{cor}

For the discrete case and $V_n = n$ we obtain:
\begin{cor} \label{c2}
{\it Let measures $\mu , \nu $ and ${\rm P}$ on  a filtered space $(\Omega ,
\mathcal{F}, \mathbf{F} )$ be such that $\mu
\stackrel{loc}{\ll } {\rm P} , \nu \stackrel{loc}{\ll } {\rm P}$.
If $T$ is a stopping time, then
\[
 \| (\mu_{T} )_a \| = {\underline\lim }_{n\to\infty \atop
\alpha\to 1-0} \left[ \sum_{k=0}^{n-1} \int_{\{ T=k \} } Y_k
(\alpha ) d{\rm P}_k + \int_{\{ n\leqslant T \} } Y_n (\alpha ) d{\rm
P}_n \right] ,
\]
 where $(\mu_{T} )_a $ is the absolutely
continuous part of $\mu_{T} $  with respect  to $\nu_{T} $. }
\end{cor}

{\bf Step 2.} {\it We compute  ${\bf E}_{\rm P} [z_T | {\cal F}_{T-} ]$, where $\mu \ll {\rm P}$}.

In the following  theorem we give a method of calculation of ${\bf E}_{\rm P} [ z_T
| {\cal F}_{T-} ]$ if we know only $z_0 $ and  the process
$z^{T-}$ interrupted at the moment $T$.

\begin{teo} \label{th5}
{\it Let $\mu\ll {\rm P}$ and $z$ be the density process of $\mu$ with respect to ${\rm P}$. If $0<
T(\omega ) <\infty $, then  for every $n\in \mathbb{N}$ denote by $k_n$ the unique natural number such that $T(\omega ) \in \left( \frac{k_n -1}{2^n } ; \frac{k_n }{2^n
} \right] $. Then (${\rm P}$-a.s.)
\begin{equation} \label{f25}
{\bf E}_{\rm P} [ z_T | {\cal F}_{T-} ] = \left\{
\begin{array}{l}
z_T , \ \omega\in \{ T=0\} \cup \{ T=\infty \} , \\
\lim_{n\to\infty } \frac{{\bf E_{\rm P}} \left[ z_{\frac{k_n}{2^n} } 1_{
\left\{  \frac{k_n -1}{2^n } < T \leqslant \frac{k_n }{2^n } \right\} }
\big | {\cal F}_{\frac{k_n -1}{2^n } } \right] }{ {\bf E_{\rm P}} \left[
1_{\left\{ \frac{k_n -1}{2^n } < T \leqslant \frac{k_n }{2^n } \right\}}
\big | {\cal F}_{\frac{k_n -1}{2^n } } \right] } ,  \   \omega \in
\{ 0< T< \infty \} .
\end{array} \right.
\end{equation}
Denote by $K_T (\omega )$ the following ${\cal F}_{T-} $-measurable
function
\[
 K_T (\omega ) = \left\{
\begin{array}{l}
1 , \ \omega\in B := \{ T=0\} \cup \{ T=\infty \} \cup \{ z_{T- }
=0 \} , \\ \lim_{n\to\infty } \frac{{\bf E_{\rm P}} \left[
z_{\frac{k_n}{2^n} } 1_{ \left\{  \frac{k_n -1}{2^n } < T \leqslant
\frac{k_n }{2^n } \right\} } \big | {\cal F}_{\frac{k_n -1}{2^n }
} \right] }{ z_{\frac{k_n -1}{2^n} } \cdot {\bf E_{\rm P}} \left[ 1_{\left\{
\frac{k_n -1}{2^n } < T \leqslant \frac{k_n }{2^n } \right\}}  \big |
{\cal F}_{\frac{k_n -1}{2^n } } \right] } , \ \omega \in
\Omega\setminus B.
\end{array} \right.
\]
Then the following equality is fulfilled
\begin{equation} \label{f26}
{\bf E}_{\rm P} [ z_T | {\cal F}_{T-} ] = z_{T-} \cdot  K_T (\omega ), \quad {\rm P}\mbox{-a.s.}
\end{equation} }
\end{teo}
Since every martingale  of class $(D)$  can be represented as a difference of two nonnegative martingales  of  class $(D)$, Theorem \ref{th5} holds for any martingale  of class $(D)$.

For  discrete time Theorem \ref{th5} is formulated as follows.
\begin{teo} \label{th6}
{\it Let $\mu\ll {\rm P}$ and $z$ be the density process  of $\mu$ with respect to ${\rm P}$. If the time-set is $\mathbb{N}$, then ($\mathrm{P}$-a.s.) }
\begin{equation} \label{f33}
{\bf E}_{\rm P} [ z_T | {\cal F}_{T-} ] = \left\{
\begin{array}{rl}
z_T , &  \omega\in \{ T=0\} \cup \{ T=\infty \} , \\ \frac{{\bf E}_{\rm P}
\left[ z_n 1_{\{ T=n \} } \big | {\cal F}_{n-1} \right] }{ {\bf E}_{\rm P}
\left[ 1_{\{ T=n \}}  \big | {\cal F}_{n-1} \right] } , & \omega \in
\{ T=n \} .
\end{array} \right.
\end{equation}
\end{teo}

It is well known that, if $T$ is a predictable  stopping time, then $ {\bf E}_{\rm P} [ z_T | {\cal F}_{T-} ] = z_{T-} .$ Example 44 of \cite[ch. V]{Del} shows that this equality is not true in the general case and moreover, $z_{T-}$ may even fail to be
integrable.  Theorem \ref{th5} gives
a simple explanation of this phenomenon (see Remark \ref{r4} after the proof of Theorem \ref{th5}).

{\bf Step 3.} {\it The general case is proved.}

We  prove Theorems \ref{th1} and \ref{th2} and Proposition \ref{p11} in Section 1. Theorems \ref{th3}-\ref{th6} are proved in Section II.

\begin{center} {\bf I. The Hahn decomposition of measures $\mu_T$ and $\nu_T$}
\end{center}

In what follows, all the equalities and the inclusions of sets are
considered up to ${\rm Q}$-null subsets.

{\bfseries\itshape Proof of Theorem \ref{th1}.} By Lemma IV.2.16 of
\cite{JS}, we have $\{ T<S\} \cap \{ h_T <\infty \} = \{ T<S\} $, ${\rm Q}$-a.s.
By the definition of $S$, $z_T \cdot z'_T >0 $ on the set $\{ T<S\} $.
Hence $\mu_T \sim {\rm Q}_T \sim \nu_T $ on the set $\{ T<S\} $.

By (\ref{e01}),  $z_T =z_{\infty }
>0$ and $z'_T =z'_{\infty } >0$ on the set  $E_1 := \{ T=S, \
T=\infty \} \cap \{ h_T < \infty \} $. Thus we have $\mu_T \sim
{\rm Q}_T \sim \nu_T $ on the set  $E_1 $.

Put $E_2 :=\{ S<T \} \cup \{ S\leqslant T , \ T<\infty \} $. Then,
by the definition of $S$ and \cite[Lemma III.3.6]{JS}, we have $z_T \cdot z'_T =0 $ on $E_2$ and hence $E_2
\subseteq \{ z_T =0\} \cup \{ z'_T =0\} $, ${\rm Q}$-a.s. Since $\mu (\{ z_T =0 \}
) = \nu (\{ z'_T =0 \} )=0,$ then $\mu_T \perp \nu_T $ on the set
$E_2 $.

By (\ref{e01}), we have $z_T \cdot z'_T =0 $ on the set
$\{ h_T =\infty \} $. Hence $\mu_T \perp \nu_T $ on this set.
Theorem \ref{th1} is proved. $\Box $

{\bfseries\itshape Proof of Theorem \ref{th2}.} {\bf 1.} Let $H, T_Y , T_M , T_z $ and
$T_{z'} $ be the stopping times of the processes $A, Y, M, z$ and $z' $
respectively. Since $z+z' =2$, then $T_z = T_{z' }$ and $T_Y
\leqslant T_z $. By the uniqueness of the Doob-Meyer decomposition, we
have $T_M \leqslant T_Y \mbox{ and } H\leqslant T_Y $. So
\begin{equation} \label{f6}
T_M \leqslant T_Y,  H\leqslant T_Y  \mbox{ and } T_Y \leqslant T_z = T_{z' }.
\end{equation}

Let $\mu =\mu_a +\mu_s , \ \nu =\nu_a +\nu_s ,$ where $\mu_a \sim
\nu_a ,\ \mu_s \perp \nu , \ \mu\perp \nu_s ,$  be the Lebesgue
 decomposition of the measures $\mu $ and $\nu $. Then $z= z_a + z_s , \
z' = z'_a + z'_s $, where $z_a , z_s  ,  z'_a , z'_s $ are the
density processes of the corresponding measures with respect to ${\rm Q}$.
Hence
\begin{equation} \label{f7}
Y=\sqrt{ (z_a + z_s )( z'_a + z'_s )} \; \mbox{ and }  \quad Y_{\infty }
=\sqrt{z_{a\infty } z'_{a\infty } } .
\end{equation}

Let $T$ be a stopping time such that $H\leqslant T$, ${\rm Q}$-a.s. Since $Y$ belongs to class $(D)$, then
\[
{\bf E}_{\rm Q} [ M_0 ] ={\bf E}_{\rm Q} [Y_T ] + {\bf E}_{\rm Q} [A_T ] = {\bf E}_{\rm Q} [Y_T ]
+ {\bf E}_{\rm Q} [A_H ] ={\bf E}_{\rm Q} [Y_H ] + {\bf E}_{\rm Q} [A_H ] .
\]
Hence
\begin{equation} \label{f8}
{\bf E}_{\rm Q} [Y_T ] = {\bf E}_{\rm Q} [Y_H ] .
\end{equation}
Since $Y$ is a supermartingale,  (\ref{f8}) yields
\begin{equation} \label{f9}
Y_U = {\bf E}_{\rm Q} [ Y_T | {\cal F}_U ] \; , \quad \forall\ H\leqslant
U\leqslant T.
\end{equation}
Putting $T\equiv\infty $  and $U=H$ in (\ref{f9}), by (\ref{f7}), we  obtain
\begin{equation} \label{f10}
\sqrt{ (z_{aH} + z_{sH} )( z'_{aH} + z'_{sH} )} = {\bf E}_{\rm Q}
[\sqrt{z_{a\infty } z'_{a\infty } } | {\cal F}_H ] \leqslant
\sqrt{z_{aH} z'_{aH} } .
\end{equation}

Let $\mu_a \not= 0$. Then (\ref{f10}) yields (${\rm Q}$-a.s.)
\begin{equation} \label{f11}
z_{sH} \cdot z'_{H} = z_{H} \cdot z'_{sH} =0.
\end{equation}
Thus for any stopping times $U$ and $T$ such that $H\leqslant U\leqslant T$ (${\rm Q}$-a.s.) we have
\begin{equation} \label{f12}
Y_{aU} = {\bf E}_{\rm Q} [ Y_{aT} | {\cal F}_U ] .
\end{equation}
Let $Z=\frac{d\mu_a }{d\nu_a } =\frac{z_a }{z'_a }$ be the density
process of the measure $\mu_a $ with respect to $\nu_a $ (we remind that
0/0 =0). Then $Z$ is a $\nu_a $-martingale of class $(D)$ and
equality (\ref{f12}) is equivalent to
\[
\sqrt{Z}_U = {\bf E}_{\nu_a } [ \sqrt{Z}_T | {\cal F}_U ] \; ,
\quad \forall\ H\leqslant U\leqslant T.
\]
Therefore $\sqrt{Z} $ and $Z$ are $\nu_a $-martingales beginning with the moment $H$. This is possible only if ($\nu_a $-a.s.)
\begin{equation} \label{f13}
Z=Z^H .
\end{equation}

By (\ref{f11}) we have  $z_a + z'_a =2$, $\nu_a $-a.s., on the set
$\mbox{\boldmath$[$} H, \infty \mbox{\boldmath$[$}$.
Hence (\ref{f13}) yields
\[
Y=\frac{1}{\sqrt{Z_H }} z_a \; \mbox{ on the set }
\mbox{\boldmath$[$} H, \infty \mbox{\boldmath$[$} \quad (\nu_a
-\mbox{a.s.}) .
\]
Therefore, for all  $t\geqslant H$ we have  ($\nu_a $-a.s.)
\[
2= (z_a )_t + (z'_a )_t = (z_a )_t \left(1+ \frac{1}{Z_H } \right)
\]
and $(z_a )_t \in (0; 2)$ does not depend on $t$. By (\ref{f11}),
the equality  $\{ 0< z_{aH} <2 \} = \{ Y_H > 0 \} $ holds ${\rm
Q}$-a.s. Hence, by (\ref{f6}), (${\rm Q}$-a.s.)
\[
H=T_Y =T_z =T_{z' } \; \mbox{ on the set  } \{ 0< z_{aH} <2 \} = \{
Y_H > 0\} .
\]

If $\omega\in \{ Y_H =0 \} $, then either $z_H =0$ or $z'_H =0$.
Hence Lemma III.3.6 of \cite{JS} yields $H=T_Y = T_z =T_{z' }$ on $\{ Y_H =0 \}$ (${\rm Q}$-a.s.). Thus $H=T_Y = T_z =T_{z' }$ ($\mathrm{Q}$-a.s.).

Let now $\mu_a =0$. Then (\ref{f9}) and (\ref{f10}) yield $Y=0$ (${\rm Q}$-a.s.) on the set
$\mbox{\boldmath$[$} H, \infty \mbox{\boldmath$[$}$. Hence
$T_Y\leqslant H$ ${\rm Q}$-a.s. Inequalities (\ref{f6}) and Lemma III.3.6 of \cite{JS} yields $H=T_Y = T_z =T_{z' }$ ($\mathrm{Q}$-a.s.).
 It is evidently that $T_Y \leqslant S$. Therefore
\[
H\leqslant S \quad \mbox{ and } \quad \{ H<S\} = \{ 0< z_H < 2, H<\infty \}
\subset \{ S=\infty \},
\]
and (\ref{f3}) holds.

It remains to prove that $H$ is the stopping time of $h$. It follows
from (\ref{f1}) and (\ref{f3}) that $h=h^H $. On the other hand,
if $h=h^T $, then (\ref{f1}), (\ref{f3}) and Lemma IV.2.16 of
\cite{JS} yield $A=A^T $ on the set $\{ T<H \} $. Then, by the
definition of $H$, ${\rm Q} (\{ T<H \} )=0$. Hence there exists the
stopping time of $h$ and it is equal to $H$.

{\bf 2.} It is an evident consequence of item 1.

{\bf 3.} It is enough to show that the set  $B=\{ 0< z_\infty <2\}$ satisfies statement 3 of the theorem.

{\it We claim that $B\in \mathcal{F}_H$.} Indeed, by  the definition of $H$, $B=\{ 0< z_H <2\}$. Thus
\[
B=B_1 \sqcup B_2, \mbox{ where } B_1 = \{ 0< z_H <2, H<\infty \} \mbox{ and } B_2 = \{ 0< z_H <2, H=\infty \}.
\]
It is enough to prove that $B_1, B_2 \in \mathcal{F}_H$. Since $B_1 = \{ 0< z_H 1_{\{ H<\infty \} } <2\}$, then, by Proposition I.1.21 of \cite{JS}, $B_1 \in \mathcal{F}_H$. By  \cite[I.1.17]{JS}, we have
\[
\left\{ \frac{1}{n} \leqslant z_k \leqslant 2-\frac{1}{n} \right\} \cap \{ k\leqslant H\} \in \mathcal{F}_H, \; \forall n,k \in \mathbb{N}.
\]
So
\[
B_2 = \bigcup_{n=1}^\infty \bigcap_{k=1}^\infty \left\{ \frac{1}{n} \leqslant z_k \leqslant 2-\frac{1}{n} ,  k\leqslant H \right\} \in \mathcal{F}_H.
\]

{\it Let us prove that $S=H_{B^c}$} (${\rm Q}$-a.s.).
By Theorem \ref{th1} (for $T=\infty$) and the definition of $H$, we
have
\[
B^c = \{ S<\infty \} \cup  \{ h_{\infty} =\infty \} =\{ S<\infty
\} \cup  \{ h_H =\infty \} =
\]
\[ \{ S<\infty \} \cup  \{ H<\infty , h_H =\infty \} \cup \{
H=\infty , h_H =\infty \}.
\]
Thus, (\ref{e01}) and (\ref{f3}) yield $\{ S<\infty \} \cup \{
H<\infty , h_H =\infty \}\subseteq \{ H=S\}$ (${\rm Q}$-a.s.). Hence, taking into account that $H\leqslant S$, we obtain $S=H_{B^c}$ (${\rm Q}$-a.s.). Theorem \ref{th2} is proved. $\Box $

{\bfseries\itshape Proof of Proposition \ref{p11}.}  We claim that $h(1)^{H_1} = h(1)$. Indeed, let $T$ be a stopping time such that $H_1 \leqslant T$ (${\rm Q}$-a.s.). Set $(n\geqslant 0)$
\[
W_n =\inf ( t :\ h(1)_t \geqslant n ) \wedge H_1 \ ,\ T_n = \left(
W_{ n}\right)_{\{ W_n < H_1 \} } \wedge T \geqslant W_n .
\]
Then $\{ T_n > W_n \} = \{ W_n = H_1 \} \cap \{ H_1 < T \} $. Thus, by the definition of $A^1 $ and the equality $\left( H_{1}\right)_{N_1 } = S_{N_1 } $ (${\rm Q}$-a.s.), $A^1_{W_n }
=A^1_{T_n } $ holds for $\omega\in \{ T_n > W_n \} $, and hence,
it holds ${\rm Q}$-a.s. Therefore, by Theorem I.3.17 of \cite{JS}, we
have
\[
{\bf E}_{\rm Q} [ h(1)_{W_n } ] = {\bf E}_{\rm Q} [ A^1_{W_n } ] = {\bf E}_{\rm Q} [
A^1_{T_n } ]  = {\bf E}_{\rm Q} [ h(1)_{T_n } ] < \infty .
\]
Since $h(1)$ is a nondecreasing process, then $ h(1)_{W_n } =
h(1)_{T_n } $. In particular,
\begin{equation} \label{f14}
 h(1)_{H_1 } = h(1)_{T} \; \mbox{ on the set } \{W_n =H_1 \} \cap \{ H_1 <
 T\} .
\end{equation}
Since $n$ is arbitrary and equality (\ref{f14}) evidently
holds on the set  $\{ H_1 = T\} $, we have
\[
h(1)_{H_1 } = h(1)_{T} \; \mbox{ on the set  } \cup_n \{W_n =H_1
\} \supseteq \{ h(1)_{H_1 } <\infty \} .
\]
If $\omega\in \{ h(1)_{H_1 } =\infty \} $, then $h(1)_{H_1 } =
h(1)_{T} =\infty $ either. Thus $h(1)_{H_1 } = h(1)_{T}$ (${\rm
Q}$-a.s.).

It remains to prove the minimality of $H_1$.  Assuming the converse we can find a stopping time $T$, $T\leqslant H_1 $,
such that $h(1)=h(1)^T$ and $ {\rm Q} ( \{ T< H_1 \} )
>0$. Then $ {\rm Q} ( \{ T< H_1 \} \cap N_1 ) >0$ since, otherwise ${\rm Q}$-a.s., $N_1 = N_1 \cap \{ H_1 =T\} \in \mathcal{F}_T$ by \cite[I.1.17]{JS} and we obtain a contradiction to the choice of
$H_1 $. Hence we have
\[
{\rm Q} ( \{ T\wedge S_n < H_1\wedge S_n \} \cap N_1 ) >0 \mbox{
for some } n.
\]
Then, by Theorem I.3.17 of \cite{JS}, we have
\[
{\bf E}_{\rm Q} [ h(1)_{T\wedge S_n } ] = {\bf E}_{\rm Q} [ A^1_{T\wedge S_n } ] <
{\bf E}_{\rm Q} [ A^1_{H_1 \wedge S_n } ] = {\bf E}_{\rm Q} [ h(1)_{H_1 \wedge S_n
} ].
\]
This contradicts to our choice of $T$. $\Box $

\begin{center}
{\bf II. Calculation of the norm of the absolutely continuous part
of $\mu_{T-}$ with respect to $\nu_{T-}$}
\end{center}

{\bf 1) The case when $T\equiv\infty$ and the time-set is
$\mathbb{N}$.}

In what follows  we need some propositions.

\begin{lem} \label{l1}
{\it Let $p$ be a Borel mapping from $(X , \mathcal{B}_X ) $ into
$(Y , \mathcal{B}_Y ) $, $\mu $  a measure on $(X ,
\mathcal{B}_X ) $ and $\alpha $ be a part of $p(\mu )$. Then
there exists the part $\mu^1 $ of $\mu $ such that $p(\mu^1 )=
\alpha $ (and $p(\mu - \mu^1 )\perp \alpha $). }
\end{lem}

\pr Set $J= \{ \gamma : \gamma $ is a part of $\mu $ such
that $p(\gamma )\perp \alpha \} $. If $J=\varnothing $, then $\mu^1
=\mu $ and $p(\mu^1 )=\alpha $. If $J\not= \varnothing $, then
any chain in $J$ is bounded. By the Zorn lemma, there exists a
maximal element which is denoted by $\mu^2 $. Evidently, this
element is unique. Set $\mu^1 =\mu - \mu^2 $. It is clear that
$\mu^1 $ is the desired part of $\mu $. $\Box $

\begin{lem} \label{l2}
{\it Let positive measures $\mu , \nu $, $\mu_0 $ and $\nu_0 $ on a measurable space $(X , \mathcal{B}_X ) $ be
such that $\mu \sim \nu $ and $\mu +\mu_0 \sim  \nu +\nu_0 $.
Then for every $\alpha \in (0;1)$ the following inequality is fulfilled}
\[
\begin{split}
 \left| \int_X \left( \frac{d(\mu + \mu_0 )}{d(\nu +\nu_0 )}
\right)^{\alpha } d(\nu +\nu_0 )\right. & \left. - \int_X \left( \frac{d\mu
}{d\nu } \right)^{\alpha } d\nu \right| \\
&
\leqslant 2 \Vert \mu
\Vert^{\alpha } \cdot \Vert \nu_0 \Vert^{1-\alpha } +2 \Vert
\mu_0 \Vert^{\alpha } \cdot \Vert \nu \Vert^{1-\alpha } +4 \Vert
\mu_0 \Vert^{\alpha } \cdot \Vert \nu_0 \Vert^{1-\alpha } .
\end{split}
\]
\end{lem}

\pr  Let us present $\mu_0 $ and $\nu_0 $ in the form
\[
\mu_0 = \mu_1 +\mu_2 , \mbox{ with } \mu_1 \ll \mu , \mu_2 \perp
\mu ,
\]
\[
 \nu_0 = \nu_1 +\nu_2 , \mbox{ with } \nu_1 \ll \nu , \nu_2 \sim
\mu_2 .
\]
Then
\[
 \left(\frac{d(\mu + \mu_0 )}{d(\nu +\nu_0
)}\right)^{\alpha }  (x) = \left(\frac{d(\mu + \mu_1 )}{d(\nu
+\nu_1 )}\right)^{\alpha }  (x) + \left(\frac{d\mu_2 }{d\nu_2
}\right)^{\alpha }  (x) , (\nu +\nu_0 )-\mbox{a.s.}
\]

Using the inequality $1\leqslant (1+x)^a \leqslant 1+ax$ (which is true for
$x\geqslant 0$ and $a\in [0;1]$), the H\"{o}lder inequality and the
fact that $\frac{d\gamma_1 }{d(\gamma_1 +\gamma_2 )}\leqslant 1 ,
(\gamma_1 +\gamma_2 )$-a.s., we obtain:
\[
 \left| \int
\left(\frac{d(\mu + \mu_0 )}{d(\nu +\nu_0 )}\right)^{\alpha }
d(\nu +\nu_0 ) - \int \left(\frac{d\mu }{d\nu }\right)^{\alpha }
d\nu \right| \leqslant
\]
\begin{equation} \label{f15}
\left| \int \left(\frac{d(\mu + \mu_1 )}{d(\nu +\nu_1
)}\right)^{\alpha } d\nu - \int \left(\frac{d\mu  }{d\nu
}\right)^{\alpha } d\nu \right|  + \int \left(\frac{d(\mu + \mu_1
)}{d(\nu +\nu_1 )}\right)^{\alpha } d\nu_1 + \int
\left(\frac{d\mu_2 }{d\nu_2 }\right)^{\alpha } d\nu_2 .
\end{equation}

Let us consider each summand separately. For the third and the second summands respectively we
have:
\begin{equation} \label{f16}
\int \left(\frac{d\mu_2 }{d\nu_2 }\right)^{\alpha } d\nu_2 \leqslant
\left( \int \frac{d\mu_2 }{d\nu_2 } d\nu_2 \right)^{\alpha } \cdot
\left( \int d\nu_2 \right)^{1-\alpha } = \Vert \mu_2 \Vert^{\alpha
} \cdot \Vert \nu_2 \Vert^{1-\alpha } \leqslant
\Vert \mu_0 \Vert^{\alpha } \cdot \Vert \nu_0 \Vert^{1-\alpha } ,
\end{equation}
\[
 \int \left(\frac{d(\mu + \mu_1 )}{d(\nu +\nu_1
)}\right)^{\alpha } d\nu_1 = \int \left(\frac{d(\mu + \mu_1
)}{d\nu_1 }\right)^{\alpha } \cdot \left(\frac{d\nu_1 }{d(\nu
+\nu_1 )}\right)^{\alpha } d\nu_1 \leqslant  \int \left(\frac{d(\mu +
\mu_1 )}{d\nu_1 }\right)^{\alpha } d\nu_1
\]
\begin{equation} \label{f17}
\leqslant \Vert \mu +\mu_1 \Vert^{\alpha } \cdot \Vert \nu_1 \Vert^{1-\alpha
} \leqslant  (\Vert \mu \Vert^{\alpha } + \Vert \mu_1 \Vert^{\alpha }
)\cdot \Vert \nu_1 \Vert^{1-\alpha } \leqslant \Vert \mu \Vert^{\alpha
} \cdot \Vert \nu_0 \Vert^{1-\alpha } + \Vert \mu_0 \Vert^{\alpha
} \cdot \Vert \nu_0 \Vert^{1-\alpha },
\end{equation}
where we used the inequality $(x+y)^{\alpha } \leqslant x^{\alpha } +
y^{\alpha } $, which holds for $x+y > 0$ and $xy\geqslant 0$, and
that $\Vert \mu + \mu_1 \Vert = \mu (X)+\mu_1 (X)$. For the
first summand in (\ref{f15}), which we denote by $I_1 $, we have
\begin{equation} \label{f18}
I_1 \leqslant \left| \int \left[ \left(\frac{d(\mu + \mu_1 )}{d(\nu
+\nu_1 )}\right)^{\alpha } - \left(\frac{d\mu  }{d(\nu +\nu_1
)}\right)^{\alpha } \right] d\nu \right|  + \left| \int \left[
\left(\frac{d\mu }{d\nu }\right)^{\alpha } - \left(\frac{d\mu
}{d(\nu +\nu_1 )}\right)^{\alpha } \right] d\nu \right|.
\end{equation}
For simplicity put $\gamma =\nu + \nu_1$. Since $\gamma \sim \mu
\sim \mu + \mu_1 \sim \nu $, for the first summand in (\ref{f18}) we
have
\[
\begin{split}
& \int \left[ \left(\frac{d(\mu + \mu_1 )}{d(\nu +\nu_1
)}\right)^{\alpha } -  \left(\frac{d\mu  }{d(\nu +\nu_1
)}\right)^{\alpha } \right] \frac{d\nu }{d(\nu +\nu_1 )}\cdot
d(\nu +\nu_1 ) \leqslant  \int \left[ \left(\frac{d(\mu + \mu_1
)}{d\gamma }\right)^{\alpha } -  \left(\frac{d\mu  }{d\gamma
}\right)^{\alpha } \right] d\gamma \\
& = \int \left(\frac{d\gamma }{d(\mu +\mu_1 )}\right)^{1- \alpha } d(\mu
+\mu_1 ) -\int \left(\frac{d\gamma }{d(\mu +\mu_1 )}\right)^{1-
\alpha } \cdot \left( 1 +\frac{d\mu_1 }{d\mu }\right)^{1- \alpha }
d\mu  \\
& \leqslant   \int  \left(\frac{d\gamma }{d(\mu
+\mu_1 )}\right)^{1- \alpha } d\mu_1
 \leqslant
\int  \left(\frac{d\gamma }{d\mu_1 }\right)^{1- \alpha } d\mu_1
\leqslant   \Vert \mu_1 \Vert^{\alpha } \cdot \Vert \gamma
\Vert^{1-\alpha } \leqslant  \Vert \mu_1 \Vert^{\alpha } \cdot (\Vert
\nu \Vert^{1-\alpha } + \Vert \nu_1 \Vert^{1-\alpha } )
\end{split}
\]
\begin{equation} \label{f19}
\leqslant \Vert \mu_0 \Vert^{\alpha } \cdot \Vert \nu \Vert^{1-\alpha } +
\Vert \mu_0 \Vert^{\alpha } \cdot \Vert \nu_0 \Vert^{1-\alpha } .
\end{equation}
For the second summund in (\ref{f18}) we have
\[
 \left| \int \left[
\left(\frac{d\mu }{d\gamma }\right)^{\alpha } \cdot
\left(\frac{d\gamma }{d\nu }\right)^{\alpha } - \left(\frac{d\mu
}{d\gamma }\right)^{\alpha } \right] d\nu \right| \leqslant \int
\left(\frac{d\mu }{d\gamma }\right)^{\alpha } \cdot \left( \left(
1+\frac{d\nu_1 }{d\nu }\right)^{\alpha } -1 \right) d\nu \leqslant
\]
\begin{equation} \label{f20}
\leqslant \alpha \int \left(\frac{d\mu }{d\gamma }\right)^{\alpha }
\cdot \frac{d\nu_1 }{d\nu } d\nu \leqslant \alpha \int \left(\frac{d\mu
}{d\nu_1 }\right)^{\alpha }d\nu_1 \leqslant \Vert \mu \Vert^{\alpha }
\cdot \Vert \nu_1 \Vert^{1-\alpha } \leqslant \Vert \mu \Vert^{\alpha }
\cdot \Vert \nu_0 \Vert^{1-\alpha } .
\end{equation}

From inequalities (\ref{f15}) - (\ref{f20}) the assertion follows.
$\Box $

The proof of the following lemma is trivial.

\begin{lem} \label{l3}
{\it Let $\mu$ be a positive measure on a measurable space $(X , \mathcal{B}_X ) $. Then for every non-negative function $f(x) \in L^1 (\mu )$, the
function $ g(\alpha ) = \int_X f^{\alpha } (x) \mu(dx) $ is
continuous on the segment $[0;1]$ and $0\leqslant g(\alpha ) \leqslant
\Vert f \Vert^{\alpha }_{L^1 } $.}
\end{lem}

The following proposition is of independent interest.

\begin{sen} \label{p2}
{\it Let $\mu$ be a positive measure on a measurable space $(X , \mathcal{B}_X ) $. Let $0\leqslant f_n (x) \to f(x), \mu$-a.s., and $\sup_n \int_X f_n (x) \mu(dx)
<\infty $. Then:

{\rm I}. The following chain of relations holds$^1$\footnotemark\footnotetext{${}^1{\overline\lim }_{n\to\infty \atop
x\to 1-0} f(x,n) = \inf_{n\in \mathbb{N} \atop \alpha \in (0,1)} \sup_{m\geqslant n \atop x\in [\alpha , 1)} f(x,n)$.}
\[
\begin{split}
 {\underline\lim }_{n\to\infty\atop\alpha\to 1-0} \int_X
f_n^{\alpha } (x) \mu(dx) & =\int_X f(x) \mu(dx) \leqslant {\underline\lim
}_{n\to\infty } \int_X f_n (x) \mu(dx) \\
& \leqslant
 {\overline {\lim
}}_{n\to\infty } \int_X f_n (x) \mu(dx) = {\overline {\lim
}}_{n\to\infty\atop\alpha\to 1-0} \int_X f_n^{\alpha } (x) \mu(dx) .
\end{split}
\]

{\rm II}. The following statements are equivalent
\begin{enumerate}
\item $\lim_{n\to\infty } \int_X f_n (x) \mu(dx) =\int_X f(x) \mu(dx) = d.$
\item $\lim_{n\to\infty\atop\alpha\to 1-0} \int_X f_n^{\alpha } (x)
\mu(dx) = d.$

Let $d_n =\int_X f_n (x) \mu(dx) \not= 0$ and $f(x) \not\equiv 0 \; (\mu
$-a.s.). Then 1 and 2 are equivalent to the following
\item
\begin{itemize}
      \item[a)] $ \lim_{n\to\infty }  d_n =d\not= 0 \; ;$
      \item[b)] $\frac{1}{d_n^{\alpha } } \int_X f_n^{\alpha } (x) \mu(dx) \to 1 $
                 uniformly in $n$ as $\alpha \uparrow 1$.
      \end{itemize}
\end{enumerate} }
\end{sen}

\pr  Let us prove the first equality in {\rm I}. By Fatou's lemma \cite[ch II, \S 6]{Sh}, $f(x)\in L^1 (\mu)$. For simplicity, set
$A=\int_X f(x) \mu(dx)$ and $B={\underline\lim }_{n\to\infty\atop\alpha\to 1-0} \int_X
f_n^{\alpha } (x) \mu(dx)$. Let $\epsilon
>0$. By Lemma \ref{l3}, we can choose $\alpha_0 $ such that
\[
 \left| \int_X f^{\alpha } (x) \mu(dx) - \int_X f(x) \mu(dx) \right| < \epsilon /2 \;
, \: \alpha\in ( \alpha_0 ; 1).
\]
Fixed $\alpha_1 \in (
\alpha_0 ; 1)$. By \cite[ch II, \S 6, Lemma 3]{Sh}, the sequence $\{ f^{\alpha_1
}_{n_1 } \}$ is uniformly integrable. Choose $n_1 $ such that $ | \int_X f^{\alpha_1
}_{n_1 } (x) \mu(dx) - \int_X f^{\alpha_1 }(x) \mu(dx) | < \epsilon /2 .$
Then $ | \int_X f^{\alpha_1 }_{n_1 } (x) \mu(dx) - \int_X f(x) \mu(dx) | <
\epsilon .$ Hence $A\geqslant B$.

Conversely, let $n_k \to\infty $ and $\alpha_k \to 1-0$ be such that $ \int_X
f^{\alpha_k }_{n_k } (x) \mu(dx) \to B$. Then, by the Lyapunov
inequality, $ \left[ \int_X f^{\alpha_k }_{n_{k+i} } (x) \mu(dx)
\right]^{\frac{1}{\alpha_k } } \leqslant \left[ \int_X f^{\alpha_{k+i}
}_{n_{k+i} } (x) \mu(dx) \right]^{\frac{1}{\alpha_{k+i} } } .$
Letting $i\to\infty $, by the Fatou lemma, we have: $ \left[ \int_X f^{\alpha_k } (x)
\mu(dx) \right]^{1/\alpha_k }\leqslant B .$ Letting $k\to\infty $, by
Lemma \ref{l3}, we obtain $A\leqslant B$.

The first inequality follows from the Fatou lemma and the second
one is trivial. Let us prove the last equality in {\rm I}. For
simplicity we shall denote the first limit by $C$ and the second one
by $D$. By Lemma \ref{l3}, $C\leqslant D$. Now we prove the inverse
inequality.

Let $n_k \to\infty $ and $\alpha_k \to 1-0$ be such that $ \int_X
f^{\alpha_k }_{n_k } (x) \mu(dx) \to D$,  as  $k\to\infty $.
Then, by the Lyapunov inequality, we have
\[
 \left[ \int_X f^{\alpha_k }_{n_{k} } (x) \mu(dx)
\right]^{\frac{1}{\alpha_k } } \leqslant \int_X f_{n_k } (x) \mu(dx) .
\]
Passing to the upper limit as $k\to\infty $, we have: $ D\leqslant
{\overline {\lim }}_{n\to\infty } \int_X f_{n_k } (x) \mu(dx) \leqslant C. $

Now we prove {\rm II}. The equivalence of 1 and 2 follows from I.

$2. \Rightarrow 3.$ Since 1 follows from 2,  the limit in
$a)$ exists and $d\not= 0$. Hence there exists the limit
$\lim_{n\to\infty\atop\alpha\to 1-0} d_n^{\alpha } = d\not= 0$.
Therefore there exists the limit of the fraction and
\[
\lim_{n\to\infty\atop\alpha\to 1-0} \frac{1}{d_n^{\alpha } }
\int_X f^{\alpha }_{n} (x) \mu(dx) = 1.
\]

Let $\epsilon >0$. Choose  $N$ and $\alpha_1 $ such that $ |
\frac{1}{d_n^{\alpha } } \int_X f^{\alpha }_{n} (x) \mu(dx) -1 | <
\epsilon \, , \, \forall n>N ,\, \forall \alpha\in (\alpha_1 ; 1)$.
By Lemma \ref{l3}, we can choose $\alpha_0 > \alpha_1 $ such that
$ | \frac{1}{d_n^{\alpha } } \int_X f^{\alpha }_{n} (x) \mu(dx) -1 | <
\epsilon \, , \, \forall n=1,\dots , N,\, \forall \alpha\in (\alpha_0 ; 1).$ The last two
inequalities prove  item $b)$.

$3. \Rightarrow 2.$ By item $a)$, we have
$\lim_{n\to\infty\atop\alpha\to 1-0} d_n^{\alpha } = d$. Hence, by
item $b)$, the limit of their product exists and
\[
\lim_{n\to\infty\atop\alpha\to 1-0} \int_X f^{\alpha }_{n} (x) \mu(dx) = d. \Box
\]

In the next proposition we find the locally absolutely continuous part $\tilde\mu $  of $\mu$ with respect to $\nu $ and justify the title ``asymptotic singular part''.

\begin{sen} \label{l4}
{\it Let $\mu $ and $ \nu $ be two probability measures on
 a filtered space $(\Omega , \mathcal{F}, \mathbf{F} )$. Denote by $\mu_t^1 $
the absolutely continuous part of $\mu_t $ with respect to $\nu_t $.
\begin{itemize}
\item[1)] Let $\mu^0_{n} $ be the unique part of $\mu$ such that $\mu^0_{n} |_{{\cal F}_n } = \mu^1_{n} $. Set
    \[
     {\tilde \mu } = \mu^0_{k} - \sum^{\infty }_{n=k} (\mu^0_{n}
    -\mu^0_{n+1} ) , \enspace \forall k\in \mathbb{N}.
    \]
    Then  $\tilde\mu $ is the locally absolutely continuous part of $\mu $
 with respect to $\nu $.
\item[2)] $\mu \stackrel{as}{\perp } \nu$  if and only if $\lim_{t\to \infty } \| \mu_t^1 \| =0$.
\item[3)] If $\alpha_0
\in (0,1)$, then
\[
\lim_{t\to\infty } \left[ \int \left(
\frac{d\mu_t }{d\nu_t } \right)^{\alpha } d\nu_t - \int \left(
\frac{d{\tilde\mu }_t }{d\nu_t } \right)^{\alpha } d\nu_t \right]
= 0\; \mbox{ uniformly on } [\alpha_0 ; 1].
\]
\end{itemize} }
\end{sen}

\pr We shall prove the proposition for discrete time only. It is clear
that the sequence $\| \mu_n^1 \| $ is not increase. Set
\[
d \equiv \lim_{n\rightarrow \infty } \Vert \mu^1_{n} \Vert = \inf
\Vert \mu^1_{n} \Vert .
\]

1) Note that $\mu^0_{n}$  is defined correctly by Lemma \ref{l1}. Then
$\mu^0_{n+1} $ is a part of $\mu^0_{n} $ and $\Vert \mu^0_{n}
\Vert = \Vert \mu^1_{n} \Vert $. Hence there exists the  limit
\[
\lim_{n\rightarrow \infty } \Vert \mu^0_{n} \Vert = \inf \Vert \mu^1_{n} \Vert = d.
\]
It is clear that ${\tilde \mu }$ is a (maybe zero) part of $\mu $ such that
\[
 \Vert {\tilde \mu } \Vert = \Vert \mu^0_{k} \Vert  -
\sum^{\infty }_{n=k} (\Vert \mu^0_{n} \Vert - \Vert \mu^0_{n+1}
\Vert ) = \lim_{n\rightarrow \infty } \Vert \mu^0_{n} \Vert =d.
\]
Since ${\tilde \mu }_n = {\tilde \mu } |_{{\cal F}_n } \ll
\mu^0_{n} |_{{\cal F}_n } = \mu^1_{n}  \ll \nu_n $ for every $n$,
we obtain that ${\tilde \mu } \stackrel{loc}{\ll } \nu .$ If $\alpha$ is a locally absolutely continuous part of $\mu $, then $\alpha $ is a part of $\mu^0_{n}$ such that $\alpha \perp (\mu^0_{n}- \mu^0_{n+1})$ for every $n\geqslant 1$. Thus, by construction, $\alpha$ is a part of $\tilde\mu $. Hence $\tilde\mu $ is the locally absolutely continuous part of $\mu $
 with respect to $\nu $.

2) follows from the equality $\| {\tilde \mu }\| =\lim_{n\rightarrow \infty } \Vert \mu^1_{n} \Vert$ in the proof of item 1).

3) Set $ {\tilde I}_n (\alpha ) = \int \left( \frac{d{\tilde \mu
}_n }{d\nu_n } \right)^{\alpha } d\nu_n .$

If $\mu_N^1  =0$ for some $N$, then $\mu_n^1 =0$ and $I_n (\alpha) ={\tilde I}_n (\alpha )$,  for every $n\geqslant N$  and $\alpha \in (0,1]$. So the statement is trivial in this case.

Assume that $\mu_n^1 \not=0$ for every natural $n$.
Let  $\nu^1_{n} $ be
the absolutely continuous part of  $\nu_n $ with respect to $\mu_n $. Let us present the measures $\mu^1_{n} $ and $\nu^1_{n} $ in the
form:
\[
 \mu^1_{n} = {\tilde \mu }_{n} + {\tilde \mu }^1_{n} +
{\tilde \mu }^2_{n} , \mbox{ with } {\tilde \mu }^1_{n} \ll
{\tilde \mu }_{n} , {\tilde \mu }^2_{n} \perp {\tilde \mu }_{n},
\]
\[
\nu^1_{n} = {\tilde \nu }^1_{n} + {\tilde \nu }^2_{n} ,
\mbox{ with } {\tilde \nu }^1_{n} \sim {\tilde \mu }_{n} ,
{\tilde \nu }^2_{n} \sim {\tilde \mu }^2_{n} .
\]

Then ${\tilde \mu }^1_{n} + {\tilde \mu }^2_{n} = \mu^1_n -
{\tilde \mu}_{n} =(\mu^0_{n} - {\tilde \mu } ) |_{{\cal F}_n }$.
Therefore
\begin{equation} \label{f21}
\lim_{n\rightarrow \infty } \Vert {\tilde \mu }^1_{n} + {\tilde
\mu }^2_{n} \Vert = \lim_{n\rightarrow \infty } \Vert \mu^0_n -
{\tilde \mu } \Vert =0.
\end{equation}

By Lemma \ref{l2} and the H\"{o}lder inequality, we obtain
\[
\begin{split}
 \vert I_n (\alpha ) - {\tilde I}_n (\alpha ) \vert & = \left| \int
\left( \frac{d\mu_n }{d\nu_n } \right)^{\alpha } d\nu_n - \int
\left( \frac{d{\tilde \mu }_n }{d\nu_n } \right)^{\alpha } d\nu_n
\right|  \\
& \leqslant \left| \int \left( \frac{d({\tilde \mu }_n + {\tilde \mu }^1_n
)}{d{\tilde \nu}^1_n } \right)^{\alpha } d{\tilde \nu }^1_n -
\int \left( \frac{d{\tilde \mu }_n }{d{\tilde \nu}^1_n }
\right)^{\alpha } d{\tilde \nu }^1_n \right| + \int \left(
\frac{d{\tilde \mu }^2_n }{d{\tilde \nu}^2_n } \right)^{\alpha }
d{\tilde \nu }^2_n  \\
& \leqslant
 2\Vert {\tilde \mu }^1_n
\Vert^{\alpha } \cdot \Vert {\tilde \nu }^1_n \Vert^{1-\alpha } +
\Vert {\tilde \mu }^2_n \Vert^{\alpha } \cdot \Vert {\tilde \nu
}^2_n \Vert^{1-\alpha } \leqslant 3\Vert {\tilde \mu }^1_n  + {\tilde
\mu }^2_n  \Vert^{\alpha } \cdot \Vert  \nu \Vert^{1-\alpha }.
\end{split}
\]
This inequality and (\ref{f21}) prove statement 2). The
proposition is proved.
$\Box $

{\bfseries\itshape Proof of Proposition \ref{p0}}. Setting $f_n (x)= \frac{d(\mu_n)_a }{d\nu_n } (x)$, the proposition immediately follows from Theorem C and Proposition  \ref{p2}(I).
$\Box $

{\bfseries\itshape Proof of Corollary \ref{c0}}. Set ${\cal G}_n =
{\cal F}_{V_n } , {\cal G}_{\infty } = {\cal F}_{T-} $,
$\mu'_n = \mu_{ V_n } $,$ \nu'_n = \nu_{ V_n } ,
\mu' = \mu_{T-} , \nu' = \nu_{T-} $. Since ${\cal F}_{T-} = \bigvee_n {\cal
F}_{ V_n } $, then the assertion follows from
Proposition \ref{p0} and the Jessen theorem.
$\Box$

{\bfseries\itshape Proof of Corollary \ref{c1}}. It is easy to see that ${\cal F}_T = \bigvee_n {\cal
F}_{T\wedge V_n } $ for every stopping time $T$. Set ${\cal G}_n =
{\cal F}_{T\wedge V_n } , {\cal G}_{\infty } = {\cal F}_T $,
$\mu'_n = \mu_{T\wedge V_n } $,$ \nu'_n = \nu_{T\wedge V_n } ,
\mu' = \mu_T , \nu = \nu_T $. Then the assertion follows from
Proposition \ref{p0} and Jessen's theorem. $\Box $

{\bf 2) Computation of ${\bf E}_{\rm P} [ z_T | {\cal F}_{T-} ]$. }

{\bfseries\itshape Proof of Theorem \ref{th5}.}
 The proof of the theorem is based on an approximation of the stopping time $T$
from below.

Let us consider the following sets
\[
 A_{0}^{n} = \{ 0=T \} \; ,
\quad A_{k}^{n} = \{ \frac{k-1}{2^n } <T \} \; ,\; k=1,\dots ,
n2^n +1 .
\]
It is clear that, if $k>0$, then $A_{k}^{n} \in {\cal F}_{T-} $ and
$A_{k}^{n} \supset A_{k+1}^{n} $. Set
\[
 B_{0}^{n} = A_{0}^{n} = \{ 0=T \} \; , \quad B_{n2^n +1}^{n} =
A_{n2^n +1}^{n} = \{ n<T \} ,
\]
\[
 B_{k}^{n} = A_{k}^{n} \setminus A_{k+1}^{n} = \{
\frac{k-1}{2^n } <T\leqslant
 \frac{k}{2^n }  \} \; , k=1,\dots , n2^n .
\]
Then $\{B_{k}^{n} \}_{k=0}^{n2^n +1}$ is a finite partition of $\Omega $ for every $n$. The proof
of the following lemma is trivial.

\begin{lem}  \label{l5}
{\it Set ${\cal F}_{0} ={\cal G}_{0}$ and denote by
${\cal G}_{n}$  the $\sigma $-algebra generated by the families
of sets $ {\cal F}_{0} , \; {\cal F}_{\frac{k-1}{2^n } }
\cap A_{k}^{n} , \; k= 1,\dots , n2^n +1. $ Then
\begin{itemize}
\item [1.] Every set $E\in {\cal G}_n $ can be uniquely represented
in the form
\begin{equation} \label{f27}
E= E_0 \sqcup E_1 \sqcup \dots \sqcup E_{n2^n +1} ,
\end{equation}
where  $E_0 \in {\cal F}_0 \cap B_0^n  , \; E_k \in {\cal
F}_{\frac{k-1}{2^n } }  \cap B_{k}^{n} ,\; k= 1,\dots , n2^n +1
$. \item [2.] ${\cal F}_{T-} = \vee_n {\cal G}_n $. $\Box$
\end{itemize} }
\end{lem}

 The restrictions of the measures $\mu $ and ${\rm
P}$ on ${\cal G}_n $ we denote by $\mu'_n $ and ${\rm P'}_n $
respectively. Decomposition (\ref{f27}) shows that
\begin{equation} \label{f28}
\mu'_n = \sum_{k=0}^{n2^n +1} \mu |_{B_k^n \cap {\cal
F}_{\frac{k-1}{2^n } }} ,\quad\quad {\rm P'}_n = \sum_{k=0}^{n2^n
+1} {\rm P} |_{B_k^n \cap {\cal F}_{\frac{k-1}{2^n } }} .
\end{equation}

We need the following lemma.
\begin{lem} \label{l6}
{\it Let $\mu $ and $\nu $ be  measures on a measurable space $(\Omega , {\cal F})$
such that $\mu\ll\nu $. Let $A\in {\cal F}$ and let ${\cal G}$ be a
$\sigma $-subalgebra of ${\cal F}$. Denote by $\mu' $ and $\nu' $
the restrictions of the measures $\mu |_{\Omega\setminus A} $ and
$\nu |_{\Omega\setminus A} $ onto the $\sigma $-algebra ${\cal G}
\cap (\Omega\setminus A)$. Then $\mu' \ll \nu' $ and ($\nu' $-a.s.)
\[
\frac{d\mu' }{d\nu' } = \frac{ {\bf E}_{\nu } [ z | {\cal G} ] -
{\bf E}_{\nu } [ z\cdot 1_A | {\cal G} ] }{ {\bf E}_{\nu } [ 1 |
{\cal G} ] - {\bf E}_{\nu } [ 1_A | {\cal G} ] } \bigg |_{
\Omega\setminus A} = \frac{ {\bf E}_{\nu } [ z\cdot
1_{\Omega\setminus A} | {\cal G} ] }{ {\bf E}_{\nu } [
1_{\Omega\setminus A} | {\cal G} ] } \bigg |_{ \Omega\setminus A} ,
\mbox{ where } z=\frac{d\mu }{d\nu } .
\] }
\end{lem}

\pr The proof we separate into two steps.

I. Define measures $P$, $Q$, and $\widetilde{\nu}$ on $(\Omega , {\cal G})$ putting
\[
P(E) =\mu (E\setminus A) ,\ Q(E) =\nu (E\setminus A) , \
{\tilde\nu } (E) =\nu (E) , \mbox{ for every } E\in {\cal G}.
\]
 Then $P\ll Q\ll {\tilde\nu }$ and
\begin{equation} \label{f29}
\frac{d P}{d Q} (\omega ) = \frac{ {\bf E}_{\nu } [ z | {\cal G} ]
- {\bf E}_{\nu } [ z\cdot 1_A | {\cal G} ] }{ {\bf E}_{\nu } [ 1 |
{\cal G} ] - {\bf E}_{\nu } [ 1_A | {\cal G} ] } (\omega ), \quad
\tilde\nu -\mbox{a.s.}
\end{equation}
Indeed, if $E\in {\cal G}$, then
\[
\begin{split}
 P(E) & =\mu (E) - \mu (E\cap A)
= \int_E {\bf E}_{\nu } [ z | {\cal G} ] d\tilde\nu  - \int_E
z\cdot 1_A d\nu \\
& =  \int_E {\bf E}_{\nu } [ z | {\cal G} ]
d\tilde\nu  - \int_E {\bf E}_{\nu } [ z\cdot 1_A | {\cal G} ]
d\tilde\nu  = \int_E \left( {\bf E}_{\nu } [ z | {\cal G} ] - {\bf
E}_{\nu } [ z\cdot 1_A | {\cal G} ]\right) d\tilde\nu .
\end{split}
\]
The analogous
calculation for $Q$ gives  $ Q(E) = \int_E \left( {\bf E}_{\nu } [ 1
| {\cal G} ] - {\bf E}_{\nu } [  1_A | {\cal G} ]\right) d\tilde\nu $.
Then equality (\ref{f29}) follows from  Lemma of \cite[ch.II, \S 7(8)]{Sh}.

II. It remains to prove that
\begin{equation} \label{f30}
\frac{d\mu' }{d\nu' } = \frac{d P}{d Q} \bigg |_{\Omega\setminus
A}, \nu'-\mbox{a.s.}
\end{equation}

If we put $X_1 =X_2 = \Omega\setminus A ,$ $ {\cal F}_1 = {\cal F}
\cap X_1 ,$ ${\cal F}_2 = {\cal G} \cap X_2 $, $Y_1 =Y_2 =\Omega
$, $ {\cal G}_1 ={\cal F} ,$ and ${\cal G}_2 ={\cal G}$, then
(\ref{f30}) follows from the next statement:

{\it Let  the following  diagram is commutative
\[
\begin{array}{ccc}
( Y_1 , {\cal G}_1 ) & \stackrel{i_2 }{\longrightarrow} & ( Y_2 ,
{\cal G}_2 ) \\ \uparrow\lefteqn{\pi_1 } &&
\uparrow\lefteqn{\pi_2} \\ ( X_1 , {\cal F}_1 ) & \stackrel{i_1
}{\longrightarrow} & ( X_2 , {\cal F}_2 )
\end{array}
\]
Assume that $X_2 \subset Y_2 ,$ ${\cal F}_2  = {\cal G}_2  \cap X_2 $ and
$\pi_2 $ is the embedding. Then for any measures $\mu ,\nu ,
\mu\ll\nu ,$ on $(X_1 , {\cal F}_1 )$ the following equality is
fulfilled ($i_1 (\nu )$-a.s.)
\begin{equation} \label{f31}
\frac{d i_1 (\mu )}{d i_1 (\nu )} = \frac{ d ( i_2 \circ \pi_1 )
(\mu )}{
 d ( i_2 \circ \pi_1 ) (\nu )} \bigg |_{X_2} .
\end{equation} }

Indeed, let $E\in {\cal F}_2 $ and $E' \in {\cal G}_2 $ be such
that $E' \cap X_2 = E$ (i.e. $\pi_2^{-1} (E' ) =E).$ By the
formula of change of variables \cite[ch.II, \S 6(8)]{Sh} and the
equality $(i_2 \circ \pi_1 ) (\nu ) = (\pi_2 \circ i_1 ) (\nu )$,
we have
\[
\begin{split}
 i_1 (\mu ) (E) & = (\pi_2 \circ i_1 ) (\mu ) (E' ) =
(i_2 \circ \pi_1 ) (\mu ) (E' ) \\
& =\int_{E' } \frac{ d ( i_2 \circ \pi_1 ) (\mu )}{ d ( i_2 \circ
\pi_1 ) (\nu )} d ( i_2 \circ \pi_1 ) (\nu ) = \int_{E} \frac{ d
( i_2 \circ \pi_1 ) (\mu )}{ d ( i_2 \circ \pi_1 ) (\nu )} (\pi_2
(x_2 )) di_1 (\nu )
\end{split}
\]
 and  (\ref{f31}) follows. The lemma is proved. $\Box $

Now we complete the proof of Theorem \ref{th5}. By Lemma \ref{l6} and (\ref{f28}), we
have
\begin{equation} \label{f32}
\frac{ d \mu'_n }{ d {\rm P'}_n } = \sum_{k=1}^{n2^n } \frac{ {\bf
E}_{\rm P} \left[ z_{\frac{k}{2^n } } 1_{B^n_k } \big | {\cal
F}_{\frac{k-1}{2^n } } \right] }{ {\bf E}_{\rm P} \left[ 1_{B^n_k} \big |
{\cal F}_{\frac{k-1}{2^n } } \right] } 1_{B^n_k } + z_0 1_{\{
T=0\} } + z_n 1_{\{ n<T\} } .
\end{equation}
By \cite[ch VII, \S 6, Theorem 1]{Sh}, $\frac{ d \mu'_n }{ d{\rm P'}_n } \to \frac{ d \mu_{T-} }{ d
{\rm P}_{T-} } $, ${\rm P}$-a.s. Thus (\ref{f25}) follows from
(\ref{f32}). The correctness of the definition of $K_T $
 follows from (\ref{f25}).

By construction, equality (\ref{f26}) is fulfilled on the set $(\Omega \setminus B)\cup \{ T=0\} \cup \{ T=\infty\}$. It is known that  $ \{ z_{T-} = 0\} \subset \left\{ \frac{d\mu_{T-} }{d{\rm P}_{T-} } =0 \right\} $, ${\rm P}$-a.s. (which is strict in the general case). Thus equality (\ref{f26}) holds also on the set $ \{ z_{T-} = 0\}$.
 The theorem is proved.
$\Box $

\begin{rem} \label{r4} It is clear that we can represent
$z_{T-} $ on $\{ T<\infty \} $ in the form
\begin{equation} \label{f34}
z_{T-} = \lim_{n\to\infty } \frac{{\bf E}_{\rm P} \left[
z_{\frac{k_n}{2^n} } 1_{ \left\{ \frac{k_n -1}{2^n } < T \leqslant
\infty \right\} } \big | {\cal F}_{\frac{k_n -1}{2^n } } \right]
}{ {\bf E}_{\rm P} \left[ 1_{\left\{ \frac{k_n -1}{2^n } < T \leqslant \infty
\right\}} \big | {\cal F}_{\frac{k_n -1}{2^n } } \right] } .
\end{equation}

If to compare (\ref{f34}) with (\ref{f25}) we can see
essential distinctions. In (\ref{f25}) the set $\left\{ \frac{k_n
-1}{2^n } < T \leqslant \frac{k_n }{2^n } \right\} $ tends to the
"point" $\{ T(\omega )=T \} $, but in (\ref{f34}) the set
$\left\{ \frac{k_n -1}{2^n } < T \leqslant \infty \right\} $ tends to
the "interval" $\{ T(\omega ) \leqslant T \} $. Hence, if the quotient
${\rm Q} (\{ t\leqslant T\} ) / {\rm P} (\{ t\leqslant T\} )$, where ${\rm Q} =
z_{\infty } {\rm P} (=\mu )$, tends to
 $\infty $ as $t\to\infty $, then we can expect that
$z_{T-} $ is not integrable. We shall demonstrate this phenomenon on Example
44 of  \cite[ch. V]{Del}.

Let $S$ be an arbitrary real function on a measurable space $(\Omega ,{\cal B})$. Set ${\cal
F}^0_t $ (respectively, ${\cal F}^0 $) is the $\sigma $-algebra
generated by $S\wedge t$, the set $\{ S\leqslant t\} \cap {\cal B}$ and
the atom $\{ S> t\} $ (respectively, $S$ and ${\cal B}$). If ${\rm
P}$ is a probability measure on $(\Omega ,{\cal F}^0 )$, we
denote by ${\cal F}_t $ (respectively ${\cal F} $) the $\sigma
$-algebra generated by the $\sigma $-algebra ${\cal F}^0_t $
(respectively ${\cal F}^0 $) and ${\rm P}$-null sets. Let $Z$ be
a nonnegative variable with ${\bf E}_{\rm P} [ Z ]=1$. Set ${\rm Q} = Z
{\rm P}$ and let $z$ be the density process of $\mathrm{Q}$ with respect to $\mathrm{P}$.
Let $F_{\rm P} (x)$ and $F_{\rm Q} (x)$ denote the distribution functions of $S$ with respect to ${\rm P}$ and ${\rm Q}$ respectively, i.e.,
\[
F_{\rm P} (x) = {\rm P} (\{ S\leqslant x\} ), \;  F_{\rm
Q} (x) = {\rm Q} (\{ S\leqslant x\} ) =\int Z 1_{\{ S\leqslant x \} } d{\rm
P}.
\]
Let us compute
$z_{S-} $ and ${\bf E}_{\rm P} [z_S | {\cal F}_{S-} ]$.

 Since $\{ t < S\} $ is an atom of ${\cal F}_t $, we obtain
\[
z_t 1_{\{ t<S \} } = \frac{ 1-F_{\rm Q} (t)}{1-F_{\rm P} (t)}
1_{\{ t<S \} } \; , \quad {\bf E}_{\rm P} \left[ 1_{\{ t+h <S \} } \big |
{\cal F}_t \right] = \frac{ 1-F_{\rm P} (t+h)}{1-F_{\rm P} (t)}
1_{\{ t<S \} } ,
\]
\[
 {\bf E}_{\rm P} \left[ z_{t+h} 1_{\{ t+h <S \} } \big | {\cal F}_t
\right] = \frac{ 1-F_{\rm Q} (t+h)}{1-F_{\rm P} (t)} 1_{\{ t<S \} } .
\]
Therefore, for $\omega \in\left\{ \frac{k_n -1}{2^n } < S \leqslant
\frac{k_n }{2^n } \right\} $, we have
\[
 z_{S-} =\lim_{n\to\infty
} \frac{1- F_{{\rm Q}} (\frac{k_n -1}{2^n } )}{ 1- F_{{\rm P}}
(\frac{k_n -1}{2^n } ) } \; , \; {\bf E}_{\rm P} [z_S | {\cal F}_{S-} ]
=\lim_{n\to\infty } \frac{F_{{\rm Q}} (\frac{k_n }{2^n } ) -
F_{{\rm Q}} (\frac{k_n -1}{2^n } )}{F_{{\rm P}} (\frac{k_n }{2^n
} ) - F_{{\rm P}} (\frac{k_n -1}{2^n } )} .
\]
In particular, let
\[
\Omega ={\mathbb{R}}_{+} ,\; {\cal B} = \{ \varnothing
,\Omega \} , \; S(\omega ) =\omega , \; d{\rm P} = e^{-\omega } d\omega
\mbox{ and } Z= S^{-2} \cdot e^S \cdot 1_{\{ S>1\} }.
\]
Then
\[
 F_{{\rm P}} (x)= 1- e^{-x} \; , \quad F_{{\rm Q}} (x)= \left(
1-\frac{1}{x} \right) 1_{\{ 1<x \} } ,
\]
and simple computations give us
\[
 z_{S-} = e^{\omega } 1_{[0;1]} +\frac{1}{\omega } e^{\omega }
1_{(1;\infty )}  , \quad {\bf E}_{\rm P} [ z_S | {\cal F}_{S-} ] = Z.
\]
 Hence $z_{S-}$ is not integrable.$\Box $
\end{rem}

{\bf 3) General case.}
{\bfseries\itshape Proof of Theorem \ref{th3}.}  Denote by $u$ and $u'$
the density processes of $\mu $ and $\nu $ with respect to ${\rm Q}$.
By Proposition \ref{p0} and  (\ref{f32}), we obtain
\[
 \| (\mu_{T-} )_a \|
= {\underline\lim }_{n\to\infty \atop \alpha\to 1-0} \left\{ \int
z_0^{\alpha } {z'}_0^{1-\alpha } 1_{\{ T=0\} } d{\rm P}_0 + \int
z_n^{\alpha } {z'}_n^{1-\alpha } 1_{\{ n<T\} } d{\rm P}_n +
\right.
\]
\begin{equation} \label{f35}
\left. \sum_{k=1}^{n2^n } \int \frac{ \left[ {\bf E}_{{\rm Q}}
\left[ u_{\frac{k}{2^n } } 1_{B_k^n } \big | {\cal
F}_{\frac{k-1}{2^n } } \right] \right]^{\alpha } \left[ {\bf
E}_{{\rm Q}} \left[ u'_{\frac{k}{2^n } } 1_{B_k^n } \big | {\cal
F}_{\frac{k-1}{2^n } } \right] \right]^{1-\alpha } }{ {\bf
E}_{{\rm Q}} \left[ 1_{B_k^n} \big | {\cal F}_{\frac{k-1}{2^n } }
\right] } 1_{B_k^n } d{\rm Q'}_n \right\} .
\end{equation}
It is enough to prove that the integrals under signs of the sums in
(\ref{f5}) and (\ref{f35}) are equal. Let $Z$ be the density
process of ${\rm Q}$ with respect to ${\rm P}$. Then $z=u\cdot Z, z'
= u' \cdot Z$. Hence, by   \cite[III.3.9]{JS}, we can assume
that ${\rm P} ={\rm Q}$. Denote the integrand  in
(\ref{f5}) and the integral (\ref{f5}) by $f$ and  $I$
respectively. By $g$ and $J$ we denote the integrand denominator and the integral in (\ref{f35}) respectively.
Then (see the diagram in the proof of Lemma \ref{l6})
\[
\begin{split}
 J& = \int_{B_k^n } \frac{f}{g} \big |_{B_k^n} d{\rm Q'}_n = \int
\frac{f}{g} \big |_{B_k^n} d{\rm Q}_{\frac{k}{2^n }} |_{B_k^n } =
\int \frac{f}{g} \cdot 1_{B_k^n} d{\rm Q}_{\frac{k}{2^n } }\\
& =
\int {\bf E}_{\rm Q} \left[ \frac{f}{g} \cdot 1_{B_k^n} \big | {\cal
F}_{\frac{k-1}{ 2^n }} \right] d{\rm Q}_{\frac{k-1}{2^n }} = \int
f d{\rm Q}_{\frac{k-1}{2^n } } = I.
\end{split}
\]
 The theorem is proved.$\Box $

\begin{center} {\bf Acknowledgements}
\end{center}

I am deeply indebted to the Referees for valuable comments and suggestions, which led to significant improvements of the manuscript.

{\large Department of Mathematics, Ben-Gurion University of the
Negev,

Beer-Sheva, P.O. 653, Israel}

{\it E-mail address}: $\quad$ saak@math.bgu.ac.il

\end{document}